%
%
%
%
\input amstex
\documentstyle{amsppt}
\NoBlackBoxes

\topmatter
\title Quadratic minima and modular forms\endtitle
\author Barry Brent\endauthor

\address Department of Mathematics, University of Minnesota, Minneapolis, MN 55455 \endaddress

\email barryb\@primenet.com\endemail


\abstract We give upper bounds
on the size of the gap between the constant term and the next
 non-zero Fourier 
coefficient
of an entire modular form of given weight 
for $\Gamma_0(2)$.   
Numerical
evidence indicates 
that a sharper bound holds for the weights $h \equiv 2 $
 $(\hskip -.05in \mod 4).$ 
We derive upper bounds 
for the minimum positive integer represented by
level two even positive-definite quadratic forms. 
Our data 
suggest that, for certain meromorphic modular forms and $p=2,3$, 
the $p$-order of the constant term
is related to the base-$p$ expansion of the 
order of the pole at infinity. 
\endabstract

\thanks The author is grateful to his advisor, Glenn Stevens, for 
many helpful conversations, and to Emma Previato, 
Rainer Schulze-Pillot, Masanobu Kaneko, and an anonymous referee for 
illuminating remarks. (He has appropriated some of their phrases.) He thanks 
Roger Frye for help with calculations,
using version 1.38.x of	PARI, \copyright 1989-1993 by C. Batut, D. Bernardi, 
H. Cohen, and 
M. Olivier, at the High Performance Computing Lab of Boston University. 
\endthanks


\endtopmatter

\document

\head 1. Introduction\endhead

Carl Ludwig Siegel showed in 
[Siegel 1969]  (English translation, [Siegel 1980])
that the constant terms of certain level one negative-weight 
modular forms $T_h$ are non-vanishing (`` \it Satz 2 \rm ''), and that this
implies an upper bound on the least positive exponent of a non-zero Fourier 
coefficient for  any level one entire modular form of weight $h$ with a 
non-zero constant term. 
Theta functions fall into this category. Their Fourier coefficients code up
representation numbers 
of quadratic forms. Consequently, for certain $h$, Siegel's result gives an 
upper bound on the
least positive integer represented by a
 positive-definite even unimodular quadratic form in $n = 2h$ variables. 
This bound is sharper than Minkowski's for large $n$. (Mallows, Odlyzko and 
Sloane have improved Siegel's bound in [Mallows, Odlyzko, and Sloane 1975].)

 John Hsia [private communication
to Glenn Stevens] 
suggested that Siegel's approach is workable for higher level forms. 
Following this hint, we 
constructed an analogue  of $T_h$ for $\Gamma_0(2)$,
which we denote as $T_{2,h}$. To prove $Satz$ 2, Siegel controlled the 
sign of the 
Fourier coefficients in the principal part of $T_h$. Following Siegel, we 
find upper bounds
 for the first positive exponent of a non-zero
Fourier coefficient occuring in the expansion at infinity 
of an entire modular form with a non-zero constant term for $\Gamma_0(2)$. 
The whole Siegel argument carries over 
for weights $h \equiv 0$ \newline $(\hskip -.07in \mod 4)$. 
It is not clear that Siegel's method forces the non-vanishing of 
the $T_{2,h}$ 
constant terms when $h \equiv 2$ $(\hskip -.07in \mod 4)$. 

In the latter case, we took two approaches. We used a simple trick
to derive a
bound on the size of the gap after a non-zero constant term 
in the case  $h \equiv 2 $ \newline  $\,(\hskip -.1in \mod 4)$ from
our $h \equiv 0 \,(\hskip -.1in \mod 4)$ result, avoiding the issue of 
the non-vanishing
of the constant term of $T_{2,h}$, but at the cost of a weaker estimate. 
Also (at the suggestion of Glenn Stevens), we searched for congruences 
that would imply the non-vanishing of the constant term of $T_{2,h}$. 
We found
numerical evidence that certain congruences dictate the 2- and 3-orders, 
not only of the
constant terms of the $T_{2,h}$, but of a wider class of meromorphic 
modular forms
of level $N \leq 3$. These congruences 
imply the non-vanishing of the constant term of $T_{2,h}$ for 
$h \equiv 2$ 
\newline $(\hskip -.1in \mod 4)$, but not for $h \equiv 0 \, 
(\hskip -.1in \mod 4).$

 Let us denote the vector space of entire modular forms of weight $h$ 
for $\Gamma_0(2)$ as $M(2,h)$.
In section 2, we prove that the second non-zero Fourier coefficient
of an 
an element of  $M(2,h)$ with non-zero constant term must have exponent
$ \leq \roman{dim} M(2,h) $ if $h \equiv 0 \,(\hskip -.1in \mod 4)$, or 
$ \leq 2\roman{dim} M(2,h) $ if $h \equiv 2 \,(\hskip -.1in \mod 4)$. 
(We will see that $\roman{dim} M(2,h) = 1+\left \lfloor \frac h4 \right \rfloor$.)
In section 3, we describe the numerical experiments which  indicate the
non-vanishing of the constant terms of $T_{2,h}$. Specifically, the experiments
suggest that if a  meromorphic
modular form for $\Gamma_0(N), 1 \leq N \leq 3$ with a normalized 
integral Fourier expansion
at infinity can be written as a quotient of two monomials in
Eisenstein series, then for $p=2,3$, the $p$-order of the constant term
is determined by the weight and the base-$p$ expansion of the pole-order. 
(We are aware of several papers in which base-$p$ 
expansions come up in analytical contexts, including discussions of the poles of 
coefficients of Bernoulli polynomials:
 [Kimura, 1988], [Adelberg I, Adelberg II, 1992], and [Adelberg III, 1996].)
In section 4, we prove
some of the congruences. In section 5, we apply the section 2 results to the
problem of level two quadratic minima. We state some conjectures in section 6.

The calculation of Fourier coefficients was usually done by formal 
manipulation of power series. When we could decompose a form into an infinite
product (for example, the form $\Delta^{-s}$), we 
applied the recursive relations of Theorem 14.8, [Apostol 1976], which is 
reproduced in section 2.   
\head 2. Bounds for gaps in the Fourier expansions of 
entire modular forms \endhead
Section 2.1 is introductory. We define several modular forms, 
some of which we will not need until section 3. In 2.2,
we compute the Fourier expansions of some 
higher level Eisenstein 
series. 
In section 2.3, we estimate the first positive exponent of a non-zero Fourier
coefficient in the expansion of an entire modular form for 
$\Gamma_0(2)$ with a non-zero constant term.

\subhead {\rm 2.1.}\enspace Some modular objects\endsubhead 
This section is a tour of the objects mentioned in the article. 
The main building blocks are 
Eisenstein series with known divisors and computable Fourier expansions. 

As usual, we denote by 
$\Gamma_0 (N)$ the congruence subgroup 

$$
\Gamma_0 (N) = \left\{\left(\matrix a & b\\c & d\endmatrix \right)\in\, 
SL(2,\bold Z): 
c \equiv 0\,(\mod N )\right\}
$$
and by $\Gamma(N)$ the subgroup

$$
\Gamma (N) = \left\{\left(\matrix a & b\\c & d\endmatrix \right)\in\, 
SL(2,\bold Z): \left(\matrix a & b\\c & d\endmatrix \right) 
\equiv \left(\matrix 1 & 
0\\0 & 1\endmatrix 
\right)\,(\mod N )\right\}.
$$
The vector space of entire modular forms of one variable in the upper half plane 
$\goth{H}$ of weight $h$  
for $\Gamma_0 (N)$ (``level $N$'') and trivial character, we denote by $M(N,h)$. 
We have an inclusion lattice satisfying:
$$ M(L,h) \subset M(N,h) \roman{ \,\,\,\,if\, and\, only\, if\, } L | N.$$
More particularly, any entire modular form for $SL(2,\bold Z)$ is also one for 
$\Gamma_0(2).$
The \it conductor \rm of $f$ is the least natural number $N$ such that
 $f \in M(N,h)$. 
The dimension 
of $M(N,h)$ is denoted by $r(N,h)$, or $r_{_h}$, or by $r$. 
We have the following formulas 
for positive even $h$. 
If $h \not \equiv 2  \,(\hskip -.1in \mod 12$), then

$$r(1,h) = \left\lfloor\frac{h}{12}\right \rfloor+1.$$
\flushpar
If $h \equiv 2 (\hskip -.1in \mod 12)$, then 

$$ r(1,h) = \left \lfloor\frac{h}{12} \right \rfloor.$$
\flushpar
For any positive even $h$, 
$$ r(2,h) = \left \lfloor\frac{h}{4} \right \rfloor+1.$$ 
(The level one formulas are standard. For example, see
[Serre 1973]. The level two formula can be derived by similar methods.)

The subspace of cusp 
forms in $M(N,h)$ is denoted by $S(N,h)$. We use standard notation for divisor 
sums:

$$										
\sigma_{\alpha} (n )  =   \sum_{0 < d |n} 	d^{\alpha} . 
$$
\flushpar
For complex $z$ \,satisfying $Im( z) >  0$, let 
$q=q(z)=e^{\,2 \pi i z}$.
  For positive even \newline $h \neq 2$, we denote the level one, weight $h$ 
Eisenstein series with 
Fourier expansion at infinity
$$
1 + \alpha_{_h}\sum_{n=1}^{\infty}\sigma_{_{h-1}}(n)\,q^n
$$
by $G_h$ or $G_h(z)$, where the numbers $\alpha_{_h}$  are given as follows. (For $h>0$, we  
follow 
[Serre 1973]; his $E_k$ are our $G_{2k}$.) The Bernoulli numbers $B_k$ are defined
by the expansion
$$
\frac{x}{e^x-1} = 1-\frac{x}{2}\, + \sum_{k=1}^{\infty}(-1)^{k+1}B_k\frac{x^{2k}}{(2k)!} \,.
$$
We set $\gamma_{_k} = (-1)^k\frac{4k}{B_k}$, and $\alpha_{_h} = \gamma_{_{h/2}}$ for $h>0$, 
while 
$\alpha_{_0} = 0$.
The first few $\alpha_{_h}$  \newline ($h \neq 2$) are given in the following table:

\newpage
\flushpar
\settabs 10 \columns
\+\cr
\+&&\it h & 0 & 2 & 4 & 6 & 8 & 10 & 12\cr
\+\cr
\+&&$\alpha_{_h}$ & 0 & -24 & 240 & -504 & 480 & -264 & $\frac{65520}{691}$\cr 
\+\cr\flushpar
The value of $\alpha_{_2}$ is included because, even though $G_{_2}$ is not a modular form, we 
will mention it in some of the observations.
We write $\Delta$ for the weight 12, level one cusp form 
with Fourier series
$$
\Delta = \sum_{n=1}^{\infty}\tau(n)\,q^{\,n} \,.
$$
and product expansion
$$
\Delta = q\prod_{n=1}^{\infty} (1-q^n)^{24} \,.
$$
Here, $\tau$ is the Ramanujan function. 
We denote the Klein modular 
invariant $G_4^3/\Delta$  by $j$, as usual. 
If $(N - 1)|24$, we essentially follow Apostol's notation ([Apostol 1989]), writing  
$$\phi_N (z )  =\Delta(Nz)/\Delta(z),$$ 
 $$\alpha=1/(N-1),$$ and 
$$\Phi_N = \phi_N^{\,\,\alpha}.$$ The $\Phi_N$ are  univalent 
meromorphic modular functions for $\Gamma_0 (N)$.
\endgraf

\endgraf
We define some weight $24$, level one cusp forms as follows. For 
positive integers $n,d$, let $$S_{n,d}=\Delta\left(\frac{n}{d}G_4^3+\left
(1-\frac{n}{d}\right)
G_6^2\right).$$
\endgraf
We introduce the level one functions $T_h$, which are elements in the construction of Siegel
described in Section 5.1. 
They are defined by the relation
$$
			T_h  = G_{12r -h + 2}\,\Delta^{-r}.
$$
Here $N=1$, so for even $h  > 2$, if 
$h \equiv 2 (\hskip -.1in \mod\,12)$, then $12r - h +2 = 0$, and otherwise 
$12r  - h +2 = 
14 - (h  \mod 12)$, where 
  $a \mod b = a-b \left \lfloor\frac ab \right \rfloor$, 
the least non-negative integer $A$ such that $A \equiv a \,(\hskip -.07in \mod b)$.
All poles of $T_h$ lie at infinity, and it has weight $2-h$.

We describe some level $N$ objects, $N=2,3$, using three special divisor sums:
$$\sigma^{\roman{odd}}(n) = \sum \Sb 0<d|n\\d\,\,\roman{odd}\endSb d,$$ 
$$\sigma_k^{\roman{alt}}(n) = 
\sum_{0<d|n } (-1)^{d}d^k,$$\flushpar and $$\sigma_{N,k}^{*}(n) =
\sum_{\Sb  0<d|n\\N \not\ \hskip -.01in | \frac nd 
\endSb } d^k.$$ \flushpar
Let $E_{\gamma,2}$ denote the unique normalized form in the one-dimensional space $M (2 , 2)$
 (\it i.e.\rm \,the leading coefficient in the Fourier expansion of the form is a $1$). 
The Fourier series is
$$
E_{\gamma,2} = 1+24\sum_{n=1}^{\infty} \sigma^{\roman{odd}}(n)\,q^{\,n}.\tag{2\,-1}
$$
\flushpar
$E_{\gamma,2}$ has a $\frac12$-order zero at points of 
$\goth{H}$ which are $\Gamma_0 (2)$\,-equivalent to  
$-\frac12 +  \frac12 i = \gamma$ (say).
The vector space $M (2,4)$ is spanned by two forms $E_{0,4}$ and $E_{\infty,4}$, 
which vanish with order one at the $\Gamma_0 (2)$\,-inequivalent zero and infinity cusps, 
respectively. They have Fourier expansions
$$
 			E_{0,4}  = 1 + 16\sum_{n=1}^{\infty}\sigma_3^{\roman{alt}}(n)\,q^{\,n}\tag{2-2}
$$
\flushpar
and
$$
			E_{\infty,4}  = \sum_{n=1}^{\infty}\sigma_{2,3}^{*}(n) \,q^{\,n}.\tag{2-3}
$$
\flushpar More generally, for $N=2,3$ and  even $k>2$, 
there is an Eisenstein series $E_{N,\infty,k}$ in $M(N,k)$ which 
vanishes at the 
infinity cusp, but not at cusps $\Gamma_0(N)$-equivalent to zero. 
(This exhausts the possibilities.) 
It has the Fourier expansion
$$
			E_{N,\infty,k}  = \sum_{n=1}^{\infty}\sigma_{N,k-1}^{*}(n) \,q^{\,n}.\tag{2-4}
$$
\flushpar
(With this notation, $E_{\infty,4} = E_{2,\infty,4}$.)
We write $$\Delta_2= E_{0,4}E_{\infty,4}.$$ The singleton family \{$\Delta_2$\} is a basis 
for the 
space $S(2,8)$. 

We construct a level two analogue of $j$ (distinct from $\phi_{_2}^{-1}$, which also plays this 
role):
$$
			j_{_2} =\, E_{\gamma,2}^2E_{\infty,4}^{-1} .
$$
The function $j_{_2}$ is analogous to $j$ because it is
 modular (weight zero) 
for 
$\Gamma_0 (2)$, holomorphic on 
the upper half plane, has  a
simple pole at infinity, generates the field of $\Gamma_0 (2)$\,-modular functions, and 
defines a 
bijection of a  $\Gamma_0 (2)$ fundamental set with \bf C\rm. 
We show all this in section 2.3.

Finally, we introduce analogues of the $T_h$. They are used in our
extension of Siegel's construction to level two. For $r  = r (2,h)$,
 $h\equiv 0 \,(\hskip -.1in \mod 4)$, we set
$$
			T_{2,h} =E_{\gamma,2}E_{0,4}E_{\infty,4}^{-r}.
$$
but if $h\equiv 2\, (\hskip -.1in \mod 4)$, we set
$$
			T_{2,h} =E_{\gamma,2}^2E_{0,4}E_{\infty,4}^{-1-r}.
$$
\flushpar

\subhead {\rm 2.2.}\enspace The Fourier expansions of the higher level Eisenstein series\endsubhead 
We will prove equations (2-1) and (for $N=2$) (2-4); equation (2-3) follows immediately. 
Our tools 
are results in [Schoeneberg 1974] . 
 The case $N=3$ of (2-4) 
can be proved the same way we handle $N=2$. This method also will give (2-2), but 
the calculations
are longer. Equation (2-2) can also be proved in the following way. 
For a non-zero modular form in $ M(2,h),$ 
the number of zeros
in a fundamental region is exactly $\frac h4$. 
([Schoeneberg 1974], Theorem 8, p.114.)
We check that the 
exponent of the first non-zero Fourier coefficent , if any, 
in the expansion of $G_4 - E_{0,4} - 256 E_{\infty,4}$ 
 exceeds $\frac h4 = 1 .$  This exponent counts
the number of zeros at $i \cdot \infty$. Hence
$$G_4 = E_{0,4} + 256 E_{\infty,4}.\tag{2-5}$$ 
 We deduce (2-2) from (2-3) and (2-5).

\subsubhead {\rm 2.2.1.}\enspace The modular form $E_{\gamma,2}$\endsubsubhead
Let $\zeta$ be the Riemann zeta function. Following Schoeneberg,
 let $G_2^*(z)$ be defined for $z \in \goth{H}$ by
$$
G_2^*(z)=2\zeta(2)+2\sum_{m \geq 1} \sum_{n \in  \bold Z } (mz+n)^{-2}.
$$

\flushpar
Then ([Schoeneberg 1974], p.63, equation (16)):
$$
G_2^*(z)=\frac{\pi^2}3 - 8\pi^2\sum_{n \geq 1}\sigma(n)\,e^{2\pi i z n} .
$$

\flushpar
(Here $\sigma$ is the usual sum of divisors.) For integers $N \geq 2$, let
$$
E(z,N) = NG_2^*(Nz) - G_2^*(z).
$$

\flushpar
[Schoeneberg 1974], p.177, gives the incorrect Fourier expansion
$$
E(z,N)=\frac{N-1}3\pi^2-8\pi^2\sum_{n \geq 1}\left(\sum_{\Sb  d|n\\d>0\\d \not \equiv 0 \mod N
\endSb }d\right)e^{2\pi inz}.
$$

\flushpar
Actually:
\proclaim{Proposition 2.1} $$E(z,N)=\frac{N-1}3\pi^2+8\pi^2\sum_{n \geq 1}\left
(\sum_{\Sb  d|n\\d>0
\\  d \not \equiv 0 \mod N\endSb }d\right)e^{2\pi inz}.$$\endproclaim 
We omit the proof. $E(z,N)$ belongs to $M(2,N)$ ([Schoeneberg 1974], pp. 177-178). 
We get (2-1) by 
setting $N=2$ in
Proposition 2.1 and noting that $r(2,2)=1$.
\subsubhead {\rm 2.2.2}\enspace Higher-weight Eisenstein series\endsubsubhead 
Let $\bold{\widehat  C}$ be the Riemann sphere.  Let $N$ and $k$ be integers with 
$N \geq 1, k \geq 3.$ Let $\bold m =\left(\smallmatrix m_1\\m_2\endsmallmatrix\right)$ and 
\bf a \rm $=\left(\smallmatrix a_1\\a_2\endsmallmatrix\right)$ be matrices
with entries in \bf Z \rm. Schoeneberg defines the
 \it inhomogenous Eisenstein series \rm \newline $G_{N,k,\bold a}:\goth{H}
\rightarrow  \bold {\widehat C}$ \rm as
$$
G_{N,k,\bold a}(z) = \sum_{\Sb \bold m\, \equiv \,\bold a \,\mod N\\
 \bold m \neq \bold 0\endSb}(m_1 z + m_2)^{-k} .
$$

\flushpar
If $N \geq 1$ and $k \geq 3$, then $G_{N,k,\bold a}$  has weight $k$ for 
$\Gamma(N)$ ([Schoeneberg 1974],
p.155, Theorem 1.) We put $$ \xi(t,N,k)=
\sum_{\Sb dt \equiv 1 \mod N\\d>0\endSb}\frac{\mu(d)}{d^k}.$$
Here $\mu$ is the M\"obius function. We should note that Schoeneberg uses the symbol 
$G^*$ in more than one way (differentiated by the subscripts) as we persist in 
following
 his notation. He introduces \it reduced Eisenstein series 
\rm  $G^*_{N,k,\bold a}$ for vectors \bf a \rm 
satisfying 
 $gcd(a_1,a_2,N)=1$, requiring that
$$
G^*_{N,k,\bold a}=\sum_{t \mod N}\xi(t,N,k)G_{N,k,t\bold a}.\tag{2-6}
$$

\flushpar
(This is equation (9), p.159 of
[Schoeneberg 1974], not Schoeneberg's original definition.)
 Schoeneberg introduces series indexed by 
level $N$ congruence \,subgroups $\Gamma_1$ of $SL(2,\bold Z)$ as follows. Let $\mu_1$ be the
 (finite) subgroup index $[\Gamma_1:\Gamma(N)]$, and let one coset decomposition of $\Gamma_1$ 
be $$ \Gamma_1 = \bigcup_{\nu=1}^{\mu_1}\Gamma(N)A_{\nu}.\tag{2-7}$$

\flushpar
Then for $gcd(a_1,a_2)=1$ he defines $G^*_{\Gamma_1,k,\bold a}$ as
$$
G^*_{\Gamma_1,k,\bold a}=\sum_{\nu=1}^{\mu_1}G^*_{N,k,(^tA_{\nu})\bold a}.\tag{2-8}
$$
\flushpar \bf Remark 2.1.\rm \,On pp. 161-162 of [Schoeneberg 1974], 
the author shows that  $G^*_{\Gamma_1,k,\bold a}$ is an
entire weight $k$ level $N$ modular form for $\Gamma_1$. He shows also (p.163)
 that, up to a multiplicative constant, there is only one
 $G^*_{\Gamma_1,k,\bold a}$ differing
 from $0$ at exactly those cusps that have the
form
$$
V\left(-\frac{a_2}{a_1}\right), \,\,\,V \in \Gamma_1.
$$

In view of (2-6) - (2-8), to calculate the Fourier expansion of  $G^*_{\Gamma_1,k,\bold a}$ 
it is sufficient to 
know the Fourier expansions of the $G_{N,k,\bold a}$. They are
 as follows \newline
([Schoeneberg 1974], p.157). We write $\zeta_{_N}$ $=e^{2\pi i/N}, \delta(x)=1$ 
if $x \in \bold Z,\, \delta(x)=0$ otherwise. Then we may write
$$
G_{N,k,\bold a}(z)=\sum_{\nu \geq 0} \alpha_{_\nu}(N,k,\bold a)e^{2 \pi iz\nu/N}\tag{2-9}
$$

\flushpar where
$$
\alpha_{_0}(N,k,\bold a)=\delta\left(\frac{a_1}N\right)
\sum_{\Sb m_2 \equiv a_2 \mod N\\\bold m \neq \bold 0\endSb}m_2^{-k},\tag{2-10}
$$

\flushpar
and, for $\nu \geq 1$,
$$
\alpha_{_\nu}(N,k,\bold a)=\frac{(-2\pi i)^k}{N^k(k-1)!}\sum_{\Sb m|\nu\\ \frac{\nu}m \equiv
\,a_1 \mod N \endSb}m^{k-1} \roman{sgn} \,m\, \zeta^{a_2m}_{_N}.\tag{2-11}
$$

\subsubhead {\rm 2.2.3}\enspace The modular form $E_{2,\infty,k}$\endsubsubhead 
Let $\bold u = \left( \smallmatrix 1\\0 \endsmallmatrix \right)$ 
and $\omega$ be the
leading Fourier coefficient in the expansion of $G^*_{\Gamma_0(N),k,\bold u}$\, ($N=2$ or $3$). 
By Remark 2.1, \newline $G^*_{\Gamma_0(N),k,\bold u}(i \infty) = 0$, so 
$$E_{N,\infty,k}=\frac1{\omega}G^*_{\Gamma_0(N),k,\bold u}\,\,\tag{2-12}$$
\flushpar
is a normalized modular form in $M(N,k)$ which vanishes at infinity but not at zero.
\proclaim{Proposition 2.2} The Fourier expansion at infinity of $E_{2,\infty,k}$ is
$$
		E_{2,\infty,k}  = \sum_{n=1}^{\infty}\sigma_{2,k-1}^{*}(n) \,q^{\,n}.
$$
\endproclaim
\flushpar \it Proof:\rm \,\,
We choose $\Gamma_1=\Gamma_0(2)$ and 
specialize (2-6)-(2-11) to 
this setting.
 The coset decomposition (2-7) is determined as follows. 
For $\Gamma_1=\Gamma_0(2), \,\,\mu_{_1}=2$, because ([Schoeneberg 1974], p.79) 
$[\Gamma_0(N):\Gamma(N)]=N \phi (N),$
 where $\phi$ is Euler's function. 
The matrices $\left(\smallmatrix 1&0\\0&1 \endsmallmatrix \right)$
 and $\left(\smallmatrix 1&1\\0&1 \endsmallmatrix \right)$ are 
inequivalent modulo $\Gamma(2)$, so we
have
$$
\Gamma_0(2) = \Gamma(2)\left(\smallmatrix 1&0\\0&1 \endsmallmatrix \right) \cup 
\Gamma(2)\left(\smallmatrix 1&1\\0&1 \endsmallmatrix \right).\tag{2-13}
$$ It is routine, so we omit the remainder of the calculation.
$\hskip .2in \boxed{}$
\subsubhead {\rm 2.2.4}\enspace The product expansion of $E_{\infty,4}$\endsubsubhead 
\proclaim{Proposition 2.3} The modular form $E_{\infty,4} \in M(2,4)$ has the following
product decomposition in the variable $q = \exp(2 \pi i z)$:
$$
E_{\infty,4}(z) = 
q\prod_{0<n\in 2\bold Z} \left(1-q^n \right)^8 \prod_{0<n\in \bold Z\backslash 2\bold Z} 
\left(1-q^n \right)^{-8}.\tag{2-14}
$$
\endproclaim

\flushpar \it Proof. \rm 
We begin by showing that for $Im(z) > 0$,  $$E_{\infty,4}(z) = 
\eta(2z)^{16}\eta(z)^{-8}.\tag{2-15}$$  For now, denote $\eta(2z)^{16}\eta(z)^{-8}$ as $F(z).$ 
The function $F$ 
is holomorphic on $\goth{H}$ because $\eta$ is non-vanishing on $\goth{H}$. $F$ has the product 
expansion
$$ F(z) = q\prod_{0<n\in 2\bold Z} \left(1-q^n \right)^8 \prod_{0<n\in \bold Z\backslash 2
\bold Z} 
\left(1-q^n \right)^{-8}.\tag{2-16}$$ This follows from the product expansion of $\eta$. It
 shows that $F$ has a simple zero at infinity. The number of zeros in 
a $\Gamma_0(2)$ fundamental 
set in $\goth{H}$ for a level $2$, weight $4$ modular form is one. If we showed that $F$ has
 weight $4$ for $\Gamma_0(2)$, it would follow that the divisors of $E_{\infty,4}$ and $F$ are both 
$1\cdot  i\, \infty \,\,.$ The expansion (2-16) shows that the Fourier series of $F$ is monic. 
So is that of $E_{\infty,4}$. Thus $F$ and $E_{\infty,4}$ would be monic modular forms with the 
same weight, 
level and 
divisor, hence identical. 
So, we only need to check the weight $4$ modularity of $F$ on a 
set of
generators for $\Gamma_0(2)$. One such set is $\left\{ T,V \right\}$, where 
$T=\left(\smallmatrix 1&1\\0&1\endsmallmatrix\right)$ and
 $V=\left(\smallmatrix -1&-1\\2&1\endsmallmatrix\right)$. ( [Apostol 1989], 
 \newline Th. 4.3.)

We calculate $F(T(z))/F(z)$ using the identity 
$$\eta(z+b)= e^{\pi ib/12}\eta(z).\tag{2-17}$$
We have 
$$F(T(z))/F(z)=F(z+1)/F(z)=\eta^{16}(2[z+1])\eta^{-8}([z+1])\eta^{-16}(2z)\eta^8(z)=$$
$$[e^{\pi i 2/12} \eta(2z)]^{16}[e^{\pi i /12}\eta(z)]^{-8}\eta^{-16}(2z)\eta^8(z)=
1=(0\,z+1)^4,$$
which is what we needed.

To check modularity for $V$, we use Dedekind's functional equation. 
This implies that $$\eta(V(z))=(-i-2iz)^{1/2}\eta(z).\tag{2-18}$$\flushpar \rm 
Equation (2-17) and Dedekind's equation also imply that 
$$\eta(2V(z))=\eta\left(\frac{-1}{2z+1}-1\right)= 
e^{-\pi i/12}\eta\left(\frac{-1}{2z+1}\right)=$$ \newline
$$e^{-\pi i/12}(-i(2z+1))^{1/2}\eta(2z+1)=
  e^{-\pi i/12}(-i(2z+1))^{1/2}e^{\pi i/12}\eta(2z)=$$ \newline 
$$(-i-2iz)^{1/2}\eta(2z).$$ That is: 
$$
\eta(2V(z))=(-i-2iz)^{1/2}\eta(2z).\tag{2-19}
$$
\flushpar
By (2-18) and (2-19), $$F(V(z)) = \eta^{16}(2V(z))\eta^{-8}(V(z))=
\frac{\left([-i-2z]^{1/2}\eta(2z)\right)^{16}}{\left([-i-2z]^{1/2}\eta(z)\right)^8}=
(2z+1)^4F(z).$$ That verifies the weight 4 modularity for $V$ and completes the proof.
$\hskip .1in \boxed{} $

\subhead {\rm 2.3.}\enspace Bounds for gaps in Fourier expansions of 
level two entire forms\endsubhead 
At the step Siegel called \it Satz 2 \rm,
his argument and our extension of it
 depend on separate, fortuitous 
sign properties of particular modular forms. These lucky accidents probably bear further study.

\subsubhead {\rm 2.3.1}\enspace Siegel's argument at level one\endsubsubhead 
Let us denote the coefficient of $q^n$ in the Fourier expansion of $f$ at infinity
as $c_n[f].$
Suppose that $f \in M(1,h)$ and $c_{_0}[f] \neq 0$.
 Siegel showed that $c_n[f] \neq 0$ for some positive $n \leq \roman{dim} M(1,h)
=r$ (say).
We sketch his argument. Siegel sets
$$W =W(f) = (G_{h -12r +12})^{-1}\Delta^{1-r}f.$$
$W(f)$ turns out to be a polynomial in $j$.

The normalized meromorphic form $T_h$  has a Fourier 
series of the form
$$T_h  = C_{h,-r}q^{-r} + ... + C_{h,0} + ... ,$$
with $C_{h,-r}  = 1$. 
	Siegel proves his \it Satz 1\rm, $c_{_0}[T_h f] = 0,$  by showing that 
$$T_h f = (2\pi i)^{-1}W(f  )\frac{dj}{dz}\,\,.$$
(Since the right member of this equation is the derivative of a polynomial in $j$, 
the constant term of its Fourier series is zero.) 

  Siegel then proves his  \it Satz 2: \rm $ C_{h,0} \neq 0.$ To illustrate his
approach, we present his argument specialized to weights
 $h \equiv 0 \, (\hskip -.07in \mod 12).$
Siegel employs the operator $\frac d{d \log \, q},$ which we will abbreviate as $D$.
At level one, for weights $ h \equiv 0 \,\,(\hskip -.1in\mod 12),$ we have 
$$ T_h = - \Delta^{1-r}Dj.$$
Also, $j\Delta = G_4^3.$ So:  $$ -T_h = \Delta^{1-r}Dj = 
D(\Delta^{1-r}j) - j D(\Delta^{1-r})
=  D(\Delta^{1-r}j) - j(1-r)\Delta^{-r}D(\Delta) = $$ \newline 
$$D(\Delta^{1-r}j) + (r-1)j\Delta^{-r}[-\frac 1r \Delta^{1+r} D(\Delta^{-r})]=
 D(\Delta^{1-r}j) + \frac {1-r}r j \Delta D( \Delta^{-r})$$ \newline 
$$ =  D(\Delta^{1-r}j) + \frac {1-r}r G_4^3 D( \Delta^{-r}) \, .$$
The term $D(\Delta^{1-r}j) $ is the derivative of a Fourier series, so it contributes nothing
to the constant term of $T_h.$ The Fourier series of $G_4^3$ has positive coefficients. 
The Fourier coefficients in the principal part of $D(\Delta^{-r})$ are negative, and it has
no constant term, so the constant term of $  G_4^3 D( \Delta^{-r}) $ is negative. 
For $r>1$ (the non-trivial case) it follows that $C_{h,0} < 0$.

	Siegel completes his argument as follows. Let the Fourier expansion of $f$  be
			$$f = A_0 + A_1 q + A_2 q^2 + ... ,$$
$A_0 \neq 0.$ Then by \it Satz 1, \rm 
$$0 = c_{_0}[T_h f ] = C_{h,0}A_0 + ... + C_{h,-r} A_r.$$ 
By hypothesis, $A_0 \neq 0,$ and by \it Satz 2, \rm 
$$	A_0 = -(C_{h ,0})^{-1}(C_{h ,-1}A_1 + ... + C_{h ,-r} A_r  ).$$
It follows that one of the $A_n  \,(n = 1, ..., r ) $ is non-zero. 

\subsubhead {\rm 2.3.2}\enspace Function theory at level two\endsubsubhead We
collect some familiar or easily verified facts. 
The point at infinity is represented as $i \cdot \infty$ and the extended 
upper half-plane as $\goth{H}^*$. The set of equivalence classes modulo
$\Gamma_0(2)$ in $\goth{H}^*$ we write as  $\goth{H}^*/\Gamma_0(2)$. This set
does have the structure of a genus zero
Riemann surface ([Schoeneberg 1974], pp. 91-93, 103). A set of representatives
for  $\goth{H}^*/\Gamma_0(2)$ is called a \it fundamental set \rm 
for $\Gamma_0(2)$, and a set $F$
in  $\goth{H}^*$ containing a fundamental set, such that distinct
$\Gamma_0(2)$-equivalent points in $F$ must lie on its boundary, is called 
a \it fundamental region \rm for $\Gamma_0(2)$. 
Let $S$ and $T$ be the linear fractional transformations
$S:z \mapsto \frac {-1}z$ and $T:z \mapsto z+1$. Let 
$$ R = \left\{z \in \goth{H} \,\,\roman{such\,\, that\,} 
|z| > 1, | \roman{Re  }\, z| < \frac 12\right\}.  $$
Let $V$ be the closure of $R \cup S(R) \cup ST(R)$ in the usual topology on $\bold{C}$,
and let $F_2 = V \cup \{ i \cdot \infty \}$. Then $F_2$ is a 
fundamental region for $\Gamma_0(2)$. It has two $\Gamma_0(2)$-inequivalent
cusps: zero and $i \cdot \infty$. The only non-cusp in $F_2$ fixed by a map in $\Gamma_0(2)$ is
$\gamma = - \frac 12 + \frac 12 i$.

Modular forms for $\Gamma_0(2)$ are not functions on $\goth{H}^*/\Gamma_0(2)$, but
the orders of their zeros and poles are well-defined.
We write $\roman{ord}_{z}(f)$ for the order of a zero or pole of 
a modular form $f$ at $z$. (This notation
supresses the dependence on the subgroup $\Gamma_1$ in  $SL(2,\bold Z)$ for which $f$ is
modular.)  
In non-trivial
(\it i.e. \rm even weight) cases, $\roman{ord}_{z}(f)$  at a point $z$ fixed by
an element of $\Gamma_0(2)$
 lies in  $\frac 12 \bold{Z},
\frac 13 \bold{Z}$, or $\bold{Z}$, depending upon whether $z$ is 
$SL(2,\bold Z)$-equivalent to $i$, to
$\rho = e^{2 \pi i/3}$, or otherwise.
(The fixed point $\gamma$ is $SL(2,\bold Z)$-equivalent to $i$.)
If $f$ and $g$ are meromorphic modular forms for a subgroup $\Gamma_1$ 
of finite index in  $SL(2,\bold Z)$, then 
$$\roman{ord}_{z}(f) + \roman{ord}_{z}(g) = \roman{ord}_{z}(fg).$$
The number of zeros in a fundamental set of a non-zero function
in $M(2,h)$ is $\frac h4$. To represent
 the divisor of a modular form for $\Gamma_0(2)$,  we choose a 
fundamental set $V_2$
and write a formal sum
$$\roman{div}(f) = \sum_{\alpha \in V_2} \roman{ord}_{\alpha}(f)[\alpha].$$
If $f$ and $g$ are meromorphic modular forms for $\Gamma_0(2)$ of equal weight such
that $\roman{div}(f) = \roman{div}(g)$, then $f = \lambda g$ for some constant $\lambda$.
We recall that $\roman{dim} M(N,h)$ is denoted as $r(N,h)$ and that the subspace
of cusp forms in $M(N,h)$ is denoted as $S(N,h)$.
\proclaim{Proposition 2.4} If $h$ is an even non-negative number, then
$ r(2,h) = \left \lfloor \frac h4 \right \rfloor + 1.$
\endproclaim

\flushpar \it Sketch of the proof. \rm First we note that multiplication by 
$\Delta_2 = E_{0,4} E_{\infty,4} \in S(2,8)$ is a vector space isomorphism between 
$M(2,h)$ and $S(2,h+8)$. Under the usual definitions (\it e.g.\rm\, 
[Ogg 1969], p. III-5), evaluation of a modular form 
at a cusp is a linear functional. Therefore the map
$$ \xi: M(2,h) \rightarrow \bold{C} \times \bold{C} $$ 
given by $$\xi(f) = (f(0), f(i\cdot \infty))$$
is linear with kernel $S(2,h)$. 
For $h \geq 4$, let $h=4n+2m$,  $m=0$ or $1$.
Since $E_{\infty,4}(0) \neq 0, E_{0,4}(i\cdot \infty) \neq 0,
E_{\infty,4}(i\cdot \infty) = 0, E_{0,4}(0) = 0$, and $E_{\gamma,2}$ vanishes at 
neither cusp, the values of
 $\xi(aE_{\gamma,2}^m E_{\infty,4}^n + bE_{\gamma,2}^mE_{0,4}^n)$ cover 
$\bold{C} \times \bold{C}$ as $a,b$ range over $\bold{C}$. Thus, $\xi$ is
surjective. Hence $\roman{dim} M(2,h) = 2 + \roman{dim} S(2,h)$. This fact 
allows an induction argument. One checks the initial cases by hand. For example,
a form in $M(2,h)$ has precisely one zero (with order $\frac 12$) at
a point $\Gamma_0(2)$-equivalent to $\gamma$ in a fundamental set. This fixes the divisor,
so $r(2,2)=1$.  $\hskip .1in \boxed{} $

Next, we show that $j_2= E_{\gamma,2}^2E_{\infty,4}^{-1}$ has properties analogous to 
those of $j$.
\proclaim{Proposition 2.5} The function $j_2 $
 is a modular function (weight zero modular form) for $\Gamma_0(2)$. 
It is holomorphic on $\goth{H}$ with a simple
pole at infinity. It defines a bijection of $\goth{H}/\Gamma_0(2)$ onto $\bold{C}$
by passage to the quotient.
\endproclaim

\flushpar \it  Proof. \rm The first two claims are obvious. 
  To establish the last claim,
let $f_{\lambda} = E_{\gamma,2}^2 - \lambda E_{\infty,4}$ for 
$\lambda \in \bold{C}$. Then $f_{\lambda} \in M(2,4)$. The sum of its
zero orders in a fundamental set is $1$. If $f_{\lambda}$ has multiple zeros in
a fundamental set, their must be exactly two of them at the equivalence class
 of $\gamma$, or exactly three at that of $\rho.\hskip .1in \boxed{}$  

\proclaim{Proposition 2.6} Let $f$ be meromorphic on $\goth{H}^*$. 
The following are equivalent. \newline (i) $f$ is a modular function for $\Gamma_0(2)$.
(ii) $f$ is a quotient of two modular forms for $\Gamma_0(2)$ of equal weight.
(iii) $f$ is a rational function of $j_2$.
\endproclaim

\flushpar \it Proof. \rm Clearly (iii) $\Rightarrow$ (ii) $\Rightarrow$ (i). For
$z \in \goth{H}^*$, let $[z]$ be the equivalence class of $z$ in $\goth{H}/\Gamma_0(2)$.
 By an abuse of the notation, we may take $f$ as in (i) as a function from 
$ \goth{H}^*/\Gamma_0(2)$ to $\bold{\widehat  C}$. The function $j_2$, 
also regarded in this fashion, is invertible. 
Let $\tilde f:\bold{\widehat C} \rightarrow \bold{\widehat C}$ satisfy 
$\tilde f = f \circ j_2^{-1}$. Then $\tilde f$ is meromorphic on $\bold{\widehat  C}$,
so it is rational. If $z \in \bold{\widehat C}$, let $u = j_2^{-1}(z) 
\in \goth{H}^*/\Gamma_0(2)$. Then $f(u) = f(j_2^{-1}(z)) = \tilde f(z)
= \tilde f(j_2(u))$. Thus $f$ is a rational function in $j_2. \hskip .1in \boxed{}$

Next, we differentiate $j_2$.
\proclaim{Proposition 2.7}
 For $z \in \goth{H}$, 
$$ \frac d{dz}\,j_2(z) = -2 \pi i E_{\gamma,2}(z) E_{0,4}(z) E_{\infty,4} (z)^{-1}. $$
\endproclaim

\it Proof: \rm It follows from the functional equation that
the derivative of a modular function (weight zero modular form) has weight two.
Therefore, both expressions represent weight two meromorphic modular forms 
for $\Gamma_2(0)$. 
The only poles of either function lie at infinity. On each side, 
the principal part of the Fourier expansion at infinity
 consists only of the term $-2\pi i q^{-1}$. Therefore the form
$$ \frac d{dz}\,j_2(z) +2 \pi i E_{\gamma,2}(z) E_{0,4}(z) E_{\infty,4} (z)^{-1} $$ is 
holomorphic, weight two. We find that it is zero in the same way that we established
equation (2-5).

 $\hskip .1in \boxed{} $

\subsubhead {\rm 2.3.3}\enspace Extension of Siegel's argument to level two\endsubsubhead 
We introduce an analogue of Siegel's $W$ map. For $h \equiv 0 \,(\hskip -.07in\mod 4)$ and
$f \in M(2,h)$, let $$W_2(f) = f E_{\infty,4}^{-h/4}.$$ For $h \equiv 2 \,(\hskip -.07in\mod 4)$,
let $$W_2(f) = f E_{\gamma,2} E_{\infty,4}^{-(h+2)/4}.$$
\proclaim{Proposition 2.8} If $h$ is positive, the restriction of $W_2$ to $M(2,h)$
 is a vector space isomorphism
onto the space of polynomials in $j_2$ of degree less than $r = r(2,h)$ 
($h \equiv 0 \,(\hskip -.07in\mod 4)$) or of degree between $1$ and $r$ inclusive
($h \equiv 2 \,(\hskip -.07in\mod 4)$).
\endproclaim

\flushpar \it Proof. \rm Suppose
$h \equiv 0 \,(\hskip -.07in\mod 4)$ and $f \in M(2,h)$. 
In view of Proposition 2.4, $$W_2(f) = fE_{\infty,4}^{1-r}\,.$$ For 
$d = 0, 1, ..., r-1$, the products $j_2^d E_{\infty,4}^{r-1}$ belong to $M(2,h)$. We
have $$W_2(j_2^d E_{\infty,4}^{r-1}) = j_2^d.$$
Let $Q$ be the subspace of $M(2,h)$ generated by the modular forms
$j_2^d E_{\infty,4}^{r-1}$, $d = 0, 1, ..., r-1$ and let $R$ be the space of
polynomials in $j_2$ of degree $\leq r-1$. $W_2$ carries $Q$ 
isomorphically onto $R$. Therefore, $\roman{dim} Q = r$. Hence
$Q = M(2,h)$. This proves the first claim.

Now let $h \equiv 2 \, (\hskip -.07in \mod 4)$. Then $$ W_2(f) = 
f E_{\gamma,2} E_{\infty,4}^{-r}.$$ For $d = 0, 1, ..., r-1$, the products
$j_2^d E_{\gamma,2} E_{\infty,4}^{r-1}$ belong to $M(2,h)$. We have
$$W_2(j_2^d E_{\gamma,2} E_{\infty,4}^{r-1}) = j_2^{d+1}.$$ The map $W_2$ carries
$E_{\gamma,2} Q$ isomorphically onto $j_2 R$. Therefore, $\roman{dim} E_{\gamma,2} Q = r.$
Hence $E_{\gamma,2} Q = M(2,h)$.  $\hskip .1in \boxed{} $

\proclaim{Proposition 2.9} For even non-negative $h$ and $f \in M(2,h)$,
the constant term in the Fourier expansion at infinity of $fT_{2,h}$ is
zero. 
\endproclaim

\flushpar \it Proof. \rm If $h \equiv 0 \,(\hskip -.07in\mod 4)$,
Then $$ W_2(f) \frac d{dz} j_2 = -f E_{\infty,4}^{1-r} 2 \pi i E_{\gamma,2} 
E_{0,4} E_{\infty,4}^{-1} 
= -2 \pi i f T_{2,h}.$$ If $h \equiv 2 \,(\hskip -.07in\mod 4)$, we
get the same result by a similar calculation. Thus, 
$fT_{2,h}$ is the derivative of a polynomial in $j_2$, so it
can be expressed in a neighborhood of infinity as the derivative 
with respect to $z$ of a power series in the variable $q = \exp (2 \pi i z)$.
This derivative is a power series in $q$ with vanishing constant term.
 $\hskip .1in \boxed{} $


\proclaim{Proposition 2.10} For positive $h \equiv 0 \, (\hskip -.07in 
\mod 4)$, the constant term in the Fourier expansion
at infinity of $T_{2,h}$ is non-zero.
\endproclaim
 
\flushpar \it Proof. \rm Let $u = 2 \pi i z = \log q$.
We retain the notation $D$ for the operator $\frac d{du}$, which has
the property that $D(q^n) = nq^n$.
Let $m_2 = j_2 - 64.$ Arguing as in the proof of (2-5), we see that
$E_{\gamma,2}^2 = E_{0,4} + 64 E_{\infty,4}$, so $m_2 = E_{0,4}E_{\infty,4}^{-1}.$ Thus 
 $$\frac d{dz} m_2 = \frac d{dz} j_2 = - 2 \pi i E_{\gamma,2} E_{0,4} E_{\infty,4}^{-1},$$
 so that $$D(m_2) = -E_{\gamma,2} E_{0,4} E_{\infty,4}^{-1}.$$ 
It follows that
$$ T_{2,h} = - E_{\infty,4}^{1-r} D(m_2).$$ Hence
$$ E_{\infty,4}^{1-r} D(m_2) = D(E_{\infty,4}^{1-r}m_2) - m_2 D(E_{\infty,4}^{1-r}) =
D(E_{\infty,4}^{1-r} m_2) - m_2(1-r)E_{\infty,4}^{-r}D(E_{\infty,4}) =$$ \newline
$$D(E_{\infty,4}^{1-r} m_2)+
(r-1)\,m_2E_{\infty,4}^{-r}[-\frac 1r  E_{\infty,4}^{1+r}D(E_{\infty,4}^{-r})]
=D(E_{\infty,4}^{1-r} m_2)+ \frac {1-r}r m_2 E_{\infty,4} D(E_{\infty,4}^{-r})$$
\newline
$$ = D(E_{\infty,4}^{1-r} m_2)+ \frac {1-r}r E_{0,4} D(E_{\infty,4}^{-r}).
$$
The term $D(E_{\infty,4}^{1-r} m_2)$ makes
no contribution to the constant term. Therefore the constant term of $T_{2,h}$
is the same as that of $\frac {r-1}r E_{0,4} D(E_{\infty,4}^{-r}) $.
We now examine the principal part of $D(E_{\infty,4}^{-r})$. 

An absolutely convergent monic power series can be written as 
an infinite product. 
The technique was used by Euler 
to prove the Pentagonal Number Theorem.
It has been codified as follows ( [Apostol 1976], Theorem 14.8):
\vskip .2in
\it  For a given set A and a given arithmetical function
f, the numbers $p_{A,f}(n)$ defined  
by the equation $$\prod_{n \in A}(1-x^n)^{-f(n)/n}\, = 1 + \sum_{n=1}^{\infty}p_{A,f}(n)x^n$$
satisfy the recursion formula $$np_{A,f}(n)\,= \sum_{k=1}^n f_A(k)p_{A,f}(n-k),   $$
where $p_{A,f}(0)=1$ and $$f_A(k) = \sum_{\Sb  d|k\\d \in A \endSb }f(d).$$

\rm

 Proposition 2.3 and Apostol's Theorem 14.8 together imply that, for fixed $s$,

$$
E_{\infty,4}^{-s}=q^{-s}\sum_{n=0}^{\infty} R(n)q^n,\tag{2-20}
$$ 

\flushpar
where $R(0)=1$ and $n>0$ implies that
$$
R(n)=\frac{8s}n\sum_{a=1}^n \sigma_1^{\roman{alt}}(a)R(n-a). \tag{2-21}
$$
 Because 
$\sigma_1^{\roman{alt}}(a)$ alternates sign, the alternation of the sign of $R(n)$ 
follows by an easy induction argument from (2-21). To be specific, 
$R(n)= U_n(-1)^n$ for some $U_n > 0$. Thus we may write
$$ E_{\infty,4}^{-r} = U_0(-1)^0 q^{-r} + U_1 (-1)^1 q^{1-r} + ... +
U_{r-1}(-1)^{r-1}q^{-1} + U_r(-1)^r + ...
$$  and hence $ D(E_{\infty,4}^{-r}) = $
$$ -rU_0 (-1)^0 q^{-r} + (1-r)U_1 (-1)^1 q^{1-r}
+ ... + (-1)U_{r-1}(-1)^{r-1}q^{-1} + 0 + ...
$$ \newline
$$ = V_{r}(-1)^{1} q^{-r} + V_{r-1}(-1)^{2} q^{1-r} + ... 
+ V_1(-1)^{r} q^{-1} + 0 + ... \,\, .
$$ for positive $V_n$.

On the other hand, the Fourier coefficient of $q^n,\,\, n \geq 0$, in the expansion of 
$E_{0,4}$
is $W_n(-1)^n$ for positive $W_n$, by (2-2). Thus the constant term of $E_{0,4} 
D(E_{\infty,4}^{-r})$
is $$\sum_{n=1}^r V_n (-1)^{r+1-n}W_n (-1)^n = (-1)^{r+1}\sum_{n=1}^r V_n W_n \,.
$$
and that of $T_{2,h}$ is the number 
$$\frac {r-1}r (-1)^{r+1}\sum_{n=1}^r V_n W_n \,.
$$
For weights $h \geq 4, r >1$. $\hskip .1in \boxed{} $

 The signs of the Fourier coefficients
are not as cooperative in the case \newline
$h \equiv 2\, (\hskip -.07in \mod 4)$, and
so far we do not have a result corresponding to Proposition 2.10 in this situation.

\proclaim{Theorem 2.1}
Suppose $f \in M(2,h)$ with Fourier expansion at infinity
$$f(z) = \sum_{n=0}^{\infty} A_n q^n, \,\,A_0 \neq 0.$$
If $h \equiv 0 \, (\hskip -.07in 
\mod 4)$, then some $A_n \neq 0, 1 \leq n \leq r(2,h)$.
If $h \equiv 2 \, (\hskip -.07in 
\mod 4)$, then some $A_n \neq 0, 1 \leq n \leq 2r(2,h).$
\endproclaim

\flushpar \it Proof. \rm First suppose that 
 $h \equiv 0 \, (\hskip -.07in \mod 4)$.
The argument tracks Siegel's in the level one case.
We still denote the coefficient of $q^n$ in the Fourier expansion of $f$ at infinity
as $c_n[f].$ 
The normalized meromorphic form $T_{2,h}$  has a Fourier 
series of the form
$$T_{2,h}  = C_{h,-r}q^{-r} + ... + C_{h,0} + ... ,$$
with $C_{h,-r}  = 1$. By Proposition 2.9, 
$$0 = c_{_0}[T_{2,h} f ] = C_{h,0}A_0 + ... + C_{h,-r} A_r.$$ 
By hypothesis, $A_0 \neq 0$. By Proposition 2.10, $C_{h ,0} \neq 0$, so 
$$	A_0 = -(C_{h ,0})^{-1}(C_{h ,-1}A_1 + ... + C_{h ,-r} A_r  ).$$
It follows that one of the $A_n  \,(n = 1, ..., r ) $ is non-zero.

Now suppose   $h \equiv 2 \, (\hskip -.07in \mod 4)$, $h = 4k+2$, 
$f \in M(2,h)$. For some monic $q$-series $F$ and some non-zero constant
$C_t$, $f=1 + C_t \, q^t F$. Let $g = f^2 \in M(2,2h)$.
Then $$g = 1 + 2C_t \, q^t F + C_t^2 q^{2t} F^2.$$
Since $2h \equiv 0 \, (\hskip -.07in \mod 4), t \leq r(2,2h) = 
1 +  \left \lfloor \frac {2h}4 \right \rfloor = 
1 +  \left \lfloor \frac {8k+4}4 \right \rfloor =
2k+2.
$ On the other hand, $r(2,h) = r(2,4k+2) =
1 +  \left \lfloor \frac {4k+2}4 \right \rfloor
=1+k$. $\hskip .1in \boxed{} $

The only obstacle to obtaining the bound $r+1$ instead of $2r$ in the second case is
the lack of a version of Proposition 2.10 for weights
 $h \equiv 2 \, (\hskip -.07in \mod 4)$. In section 3, we present experimental
evidence for, among other things, an extended Proposition 2.10.

While it is possible that the level two result extends to the other levels 
$N$ at which
$\Gamma_0(N)$ has genus zero
($N = 1, \, ..., \, 10,\,12,\,13,\,16,\,25$),  the question has been raised
(by Glenn Stevens) whether, because of the absence of an 
analogue for $j$, higher genus
is an obstruction to this sort of argument.

\head 3. Observations \endhead
The divisor of a meromorphic modular form $f$, normalized so that the leading Fourier 
coefficient is $1$, determines the Fourier expansions of $f$ , 
because the divisor determines $f$. 
This suggests the problem of
finding effective rules governing the map from divisors to Fourier series. 
Some results in this 
direction are known. For example,
Fourier expansions of Eisenstein series with prescribed behavior at the cusps 
are stated in [Schoeneberg, 1974].
 
Here we study rules by which the divisor governs congruences for the 
Fourier expansion.
The theory of congruences
among holomorphic modular forms is significant in number theory, so 
it is natural to scrutinize 
any new congruences among modular forms. Regularities 
among the constant terms suggest an empirical basis for such a
theory in the meromorphic setting. 

In sections 3.1 and 3.2, we discuss three rules (for conductors $N=1,2,3$) 
governing the constant term of
the Fourier expansion at infinity.
We describe numerical evidence for congruences obeyed by certain 
meromorphic modular forms.
The congruences relate geometric and arithmetic 
data: the divisor, and the 2-order or 3-order of the constant terms. 
This connection is expressed in 
terms of the weight and the sum of the digits in the base two or base three 
 expansion of the pole order. 

These rules are described for modular forms of level $N \leq 3$. They do not 
apply to
all the objects we surveyed,  
and we don't know how to sort the deviant from non-deviant forms, except by 
inspection.
The deviations are systematic in the sense that the constant terms
at a given level still obey simple rules. 
We can also manufacture linear combinations of non-deviant forms which depart from
the congruence rules in a stronger sense: the 2-order and the 3-order of the constant
terms are arbitrary. This means that
the constant terms of some of the deviant forms 
are controlled by invariants of the divisor other than the weight and the
order of the pole at infinity. 

In our surveys, a meromorphic modular form $f$ which obeys the congruences
always has a normalized rational 
Fourier expansion and a pole at infinity. 
The $T_h$ and  $T_{2,h}$  were the first examples. 
We looked for other instances of this behavior and found it exhibited by 
some standard objects. 
We then conducted a more or less systematic survey  of similar objects.

\rm We describe two sets of data. The first survey suggests rules 
regarding the the 2-order or 
3-order of constant terms of a family of level $N$ objects,
$1 \leq N \leq 3$.
The second survey looks at negative powers of the functions $E_{N,\infty,k}$, $N=2,3$. 
 These examples form
 families of their own, and within these families, the behavior of the 
constant term is again 
predictable. 

Congruences for constant terms seem to have implications for the whole Fourier
expansion of related meromorphic forms.
In section 3.3, we report observations on the Fourier expansion of $j$ that support
this idea.
  
\subhead {\rm 3.1.}\enspace Observations on the constant terms: first survey\endsubhead 
We list several thousand forms obeying rules governing their constant terms.
	Let $d_b(n)$ be the sum of the digits in 
the base\,-$b$  expansion of the positive integer 
 $n$ and  
$c_n[f\,]$ be the coefficient $c_n$ in the Fourier series 
$$
f  = \sum_n c_n\,q^{\,n}. 
$$\flushpar
Let $p$ be prime. If an integer $n$  can be factored as $n = p^a m$, $(p,m)=1$, then we 
write:
$$
		 	\roman{ord}_p(n ) = a.
$$\flushpar
In addition we write $\roman{ord}_p(0) = \infty$. If a rational number $x$ can be 
written $\frac nd$ as 
a quotient of
integers, we set $\roman{ord}_p(x) = \roman{ord}_p(n)-\roman{ord}_p(d)$.\endgraf
We write $C_{_2}$ for the set of level two meromorphic modular forms $f$ of any weight with 
rational Fourier expansion at infinity, leading coefficient 1, and a pole at infinity of order 
$s=s(f)$ $ > 0$ such that 
$$
	\roman{ord}_2(c_{_0}[f\,]) = 3d_2(s)\,.
$$
The set of level three meromorphic modular forms $f$ of any weight with 
rational Fourier expansion at infinity, leading coefficient 1, and a pole at infinity of order 
$s=s(f)$ $> 0$ such that 
$$
		\roman{ord}_3(c_{_0}[f\,]) = d_3(s)\,.
$$
will be denoted $C_3$.

For a function $f$ with a pole of order $s=s(f)$  at infinity, let $\beta = d_2(s)$ and 
$\gamma = d_3(s)$.  
Membership in $C_{_2}$ is a congruence relation, since 

$$
			\roman{ord}_2(n) = a \Leftrightarrow   n \equiv 2^a\ (\hskip -.1in\mod\ 2^{a+1}),
$$
\flushpar
but membership in $C_3$ means a choice of two congruences:
$$
			\roman{ord}_3(n) = a \Leftrightarrow   n \equiv \pm 3^a\ (\hskip -.1in\mod\ 3^{a+1}).
$$
\flushpar  
We define two subsets of $C_3$, the members of which make this choice systematically: 
$$
D_3  = \left\{\ f \in  C_3  : c_{_0}[f\,] \equiv (-1)^s 3^{\gamma}\  
(\hskip -.1in\mod\ 3^{\gamma +1})\ \right\}
$$ and
$$ E_3  = \left\{\ f \in  C_3  : c_{_0}[f\,] \equiv  3^{\gamma}\  
(\hskip -.1in\mod\ 3^{\gamma +1})\ \right\}.
$$
If $f$ is a meromorphic modular form, let $w = w(f)$ be the weight of $f$. 
As above, let $s=s(f)$ be the order of the pole of $f$ at infinity.
Finally, we will write $L=L(f)$ for the largest digit in the base\,-3
expansion of $s(f)$.

In this survey, the constant terms of the meromorphic forms we studied have 
three modes of behavior, depending upon the conductor. 
\roster\item[1]
The meromorphic forms $f$ for $SL(2,\bold Z)$ (conductor one forms) obey 
the following rule. 
\itemitem(a)
If \,$w \equiv 0 \,(\hskip -.1in \mod\,4)$, then $f \in C_2$. 
\itemitem(b)
If \,$w \equiv 2 \,(\hskip -.1in \mod\,4)$, then  $2^{4\beta}|c_{_0}[f]$. 
\itemitem(c)
If $w \equiv 0 \,(\hskip -.1in \mod\,3)$, then $f \in D_3$. 
\itemitem(d)
If $w \equiv 1 \,\,(\hskip -.1in \mod\,3)$ and $L=1$, then $f \in E_3$.
\itemitem(e)
If $w \equiv 1 \,(\hskip -.1in \mod\,3)$ and $L=2$, then
 $3^{\gamma+1}|c_{_0}[f].$
\itemitem(f)
If $w \equiv 2 \,(mod\,3)$ , then $3^{\gamma+1}|c_{_0}[f]. $

\vskip .2in

\item
Forms with conductor two  obey (a)\,-(b), but not (in general) (c)\,-(f).

\vskip .2in
\item
Forms with conductor three obey (c)\,-(f), but not (in general) (a)\,-(b).

\endroster

\subsubhead {\rm 3.2.1.}\enspace Conductor one 
\endsubsubhead
What follows is a list of objects obeying rule (1) above.  
(The function $G_2$ isn't modular in the ordinary sense, 
but we assigned it weight $2$ to see what would happen.) 
\vskip .1in
\hskip .7in $\Delta^{-a}, \,\,\,\,\,\,\,\,\,\,\,\,\,\,\,\,\,\,1 \leq a \leq 140,$
\vskip .1in
\hskip .7in $j^a, \,\,\,\,\,\,\,\,\,\,\,\,\,\,\,\,\,\,\,\,\,\,\,\,1 \leq a \leq 50,$
\vskip .1in
\hskip .7in $j\Delta^{-a}, \,\,\,\,\,\,\,\,\,\,\,\,\,\,\,1 \leq a \leq 100,$
\vskip .1in
\hskip .7in $j^a\Delta^{-b}, \,\,\,\,\,\,\,\,\,\,\,\,\,1 \leq a,b \leq 50,$
\vskip .1in
\hskip .7in $G_6^{\,\,a}\Delta^{-b}, \,\,\,\,\,\,\,\,\,\,1\leq a,b \leq 50$.
\flushpar
(If we set $a=1$, these are the functions $T_h, h \equiv 8 \, (\hskip -.07in \mod 12), 
8 \leq h \leq 596.$)
\vskip .1in
\hskip .7in $G_4^{\,\,a}G_6^b\Delta^{-c}, \,\,\,1\leq a,c \leq 50, \,\, 0 \leq b \leq 11$,
\vskip .1in
\hskip .7in $G_{10}^{\,\,a}\,\Delta^{-b}, \,\,\,\,\,\,\,\,\,1\leq a,b \leq 50$,
\vskip .1in
\hskip .7in $G_{14}^{\,\,a}\,\Delta^{-b}, \,\,\,\,\,\,\,\,\,1\leq a,b \leq 50$,
\vskip .1in
\hskip .7in $G_{2a}\,\Delta^{-b}, \,\,\,\,\,\,\,\,\,1\leq a \leq 7, 1\leq b \leq 140$.
\flushpar
(If we set $a=2$, these are the functions $T_h, h \equiv 10 \, (\hskip -.07in \mod 12), 
10 \leq h \leq 1678,$ and if we set $a=4$, they
 are the functions $T_h, h \equiv 6 \, (\hskip -.07in \mod 12), 
6 \leq h \leq 1674.$) 

\vskip .1in
\hskip .7in $G_{2a}\,\Delta^{-b}, \,\,\,\,\,\,\,\,\,8\leq a \leq 24, 1\leq b \leq 50$,
\vskip .1in
\hskip .7in $G_{2a}^{-1}\,\Delta^{-b}, \,\,\,\,\,\,\,\,1\leq a \leq 18, 1\leq b \leq 50$,
\vskip .1in
\hskip .7in $S_{n,d}^{\,-a}$,\,\,\,\,\,\,\,\,\,\,\,\,\,\,\,\,\,
\,\,\,\,$1 \leq a \leq 50,1 \leq d \leq 4, 1 \leq n \leq d,$
\vskip .1in
\hskip .7in $G_{10}^{\,\,a}S_{1,2}^{\,-b}, \,\,\,\,\,\,\,\,\,\,1\leq a,b \leq 50$,
\vskip .1in
\hskip .7in $G_{14}^{\,\,a}S_{1,2}^{\,-b}, \,\,\,\,\,\,\,\,\,\,1\leq a,b \leq 50,$
\vskip .1in
\hskip .7in $G_{4}^{\,\,a}G_{6}^{\,\,b}S_{1,2}^{\,-c}, 
\,\,1\leq a,c \leq 50, \,1 \leq b \leq 5$.
\vskip .1in
In an earlier survey, we found that $\Delta^{-s} \in C_2 \cap C_3$ for $1 \leq s \leq 3525$.
We also found that $j^s \in  C_2 \cap C_3$ for $1 \leq s \leq 200$, and 
that $j^k\Delta^{-m} \in C_2 \cap C_3$ for 
$1 \leq k, m \leq 100$.
 The 
computing power we exploited at the time (with Roger Frye's assistance) was not 
available when we were conducting the experiments described here, so we do not 
have data on membership in $D_3$ for the additional functions.

\subsubhead {\rm 3.2.2.}\enspace Conductor two 
\endsubsubhead
This is a list of objects obeying rule (2) above. The first two items are the
first few
functions $T_{2,h}$. Evidently, rule (2) does not force the vanishing of 
the constant terms of
the functions $E_{\gamma,2}E_{0,4}E_{\infty,4}^{-a} = T_{2,h}$ for 
 $h \equiv 0$ \newline $(\hskip -.07in \mod 4)$, but would imply a level 2 \it Satz 2 \rm for 
$h \equiv 2 \, (\hskip -.07in \mod 4)$ if it held for all 
$E_{\gamma,2}^{\,2}E_{0,4}E_{\infty,4}^{-a}, a$ positive.
\vskip .1in
\hskip .7in $E_{\gamma,2}E_{0,4}E_{\infty,4}^{-a},\,\,\,\,\,\,1 \leq a \leq 100$,
\vskip .1in
\hskip .7in $E_{\gamma,2}^{\,2}E_{0,4}E_{\infty,4}^{-a},\,\,\,\,\,1 \leq a \leq 100$,
\vskip .1in
\hskip .7in $j_2^{\,a}, \hskip .9in 1 \leq a \leq 100$, 
\vskip .1in
\hskip .7in $\phi_2^{\,-a},\,\,\,\,\,\,\,\,\,\,\,\,\,\,\,\,\,\,\,\,\,\,\,\,\,1 \leq a \leq 100$,
\vskip .1in
\hskip .7in $G_{2a}E_{\infty,4}^{-b}, \,\,\,\,\,\,\,\,\,\,\,\,0 \leq a \leq 24, 1 \leq b \leq 50$,
\vskip .1in
\hskip .7in $G_{2a}^{-1}E_{\infty,4}^{-b}, \,\,\,\,\,\,\,\,\,\,\,1 \leq a \leq 11, 1 \leq b \leq 50$,
\vskip .1in
\hskip .7in $G_4^{\,\,a}E_{\infty,4}^{-b},\,\,\,\,\,\,\,\,\,\,\,\,\,1 \leq a,b \leq 50$,
\vskip .1in
\hskip .7in $G_6^{\,\,a}E_{\infty,4}^{-b},\,\,\,\,\,\,\,\,\,\,\,\,\,1 \leq a,b \leq 50$,
\vskip .1in
\hskip .7in  $G_4^{\,\,a}G_6E_{\infty,4}^{-b},\,\,\,\,\,\,1 \leq a,b \leq 50$,
\vskip .1in
\hskip .7in $G_{10}^{\,\,a}E_{\infty,4}^{-b},\,\,\,\,\,\,\,\,\,\,\,\,\,1 \leq a,b \leq 50$,
\vskip .1in
\hskip .7in $E_{\gamma,2}^{\,\,a}E_{\infty,4}^{-b},\,\,\,\,\,\,\,\,\,\,\,\,\,\,1 \leq a,b \leq 50$,
\vskip .1in
\hskip .7in $E_{\gamma,2}^{\,\,a}\Delta^{-b},\,\,\,\,\,\,\,\,\,\,\,\,\,1 \leq a,b \leq 50$,
\vskip .1in
\hskip .7in $E_{0,4}^{\,\,a}E_{\infty,4}^{-b},\,\,\,\,\,\,\,\,\,\,\,\,\,\,1 \leq a,b \leq 50$,
\vskip .1in
\hskip .7in $\Delta_2^{\,\,-a}, \,\,\,\,\,\,\,\,\,\,\,\,\,\,\,\,\,\,\,\,\,1 \leq a \leq 100$.
\subsubhead {\rm 3.2.3.}\enspace Conductor three 
\endsubsubhead
This is a brief list of objects obeying rule (3) above. It should be noted that 
$\phi_3^{-1}$ has a 
double pole at infinity. More objects obeying (3) are listed in the next section.
\vskip .1in
\hskip .7in $\phi_3^{\,-a},\,\,\,\,\,\,\,\,\,\,\,\,\,\,\,\,\,\,\,\,\,\,\,\,\,\,\,1 \leq a 
\leq 100$,
\vskip .1in
\hskip .7in $G_{10}^{\,\,-a}\phi_3^{-b},  \,\,\,\,\,\,\,\,\,\,\,\,1 \leq a,b \leq 50$,
\vskip .1in
\hskip .7in $\Phi_3^{-a}, \hskip .5in \,\,1 \leq a \leq 50$,
\vskip .1in
\hskip .7in $G_{2a}\Phi_3^{-b},\hskip .28in \,\,1 \leq a \leq 24, 1 \leq b \leq 50$,
\vskip .1in
\hskip .7in $G_4^{\,\,a}\Phi_3^{-b},\hskip .28in  \,\,1 \leq a,b \leq 50$,
\vskip .1in
\hskip .7in $G_{10}^{\,\,a}\Phi_3^{-b},\hskip .28in  \,\,1 \leq a,b \leq 50$.
\vskip .1in

\subhead {\rm 3.2.}\enspace Second survey, with deviations from rules (1)-(3)\endsubhead 
Given a pair of objects of the same conductor, pole order and weight
 obeying rules (1) - (3),
 one can find a linear combination which violates
the rules. For example, let $r(N,k)>1$ and 
let $f,g$ be distinct normalized forms in $M(N,k), 1 \leq N \leq 3$. 
Let $s$ be a positive integer.
Then $\phi=f\Delta^{-s},\gamma=g\Delta^{-s}$ are normalized, and they have 
equal pole order, weight and conductor. Suppose they are subject to one of the above 
rules which dictates for $p=2$ or $3$ that 
$\roman{ord}_p(c_{_0}[\phi])=\roman{ord}_p(c_{_0}[\gamma\,])=\epsilon$ (say). 
Further, let the constant terms of $\phi$ and 
$\gamma$ be $p^{\epsilon}\frac ab, p^{\epsilon}\frac cd$ 
(with none of $a,b,c,d$ divisible by $p$).
Let $x=(p^{\sigma}-ad)/(bc-ad)$ for an arbitrary number $\sigma \neq 0$. 
Then the meromorphic modular form
$\zeta=(1-x)\phi+x\gamma$ is also normalized with the same weight and pole order as
 $\phi$ and $\gamma$. It may have lower conductor if $N \neq 1$, but whichever rule dictated
the values of 
$\roman{ord}_p(c_{_0}[\phi]), \,\roman{ord}_p(c_{_0}[\gamma\,])$ is also part of rule (1). 
Yet it fails, because 
$c_{_0}[\zeta]=p^{\epsilon+\sigma}/bd$. 
This shows that
features of the divisor other than the weight and the order 
of the pole at infinity influence the 
arithmetic of the constant term.

This fact led us to search for other deviants.
We found systematic deviations from rules (1) - (3),
but for these examples, the 2- and 3-orders of the
constant terms were still determined by the weight and the 
order of the pole at infinity.

The following functions obey rule (2):
\vskip .1in
\hskip .7in $E_{\infty,4}^{-a}, \,\,\,\,1\leq a \leq 51$, 
\vskip .1in
\hskip .7in $E_{2,\infty,k}^{-a}, \,\,\,\,1\leq a \leq 51,6 \leq k \leq 22,\, 
k \equiv 2\, (\hskip-.1in\mod 4)$, 
\vskip .1in
\hskip .7in $	E_{2,\infty,k}^{-a}, \,\,\,\,2\leq a \leq 50, a \,\roman{even},\, 
8 \leq k \leq 24, k\equiv 0\, (\hskip-.1in\mod 4)$, 

\vskip .1in \flushpar 
and the following functions obey rule (3):
\vskip .1in
\hskip .7in $E_{3,\infty,6}^{-a}, \,\,\,\, 1 \leq a \leq 98$,
\vskip .1in
\hskip .7in $E_{3,\infty,k}^{-a}, \,\,\,\,3 \leq a \leq 48,
a \equiv 0 \,(\hskip -.1in\mod 3), k \equiv 0 \,(\hskip -.1in\mod 6), \,12 \leq k \leq 24$, 
\vskip .1in
\hskip .7in $E_{3,\infty,k}^{-a}, \,\,\,\,2 \leq a \leq 98,
a \equiv 0$ or $2 \,(\hskip -.1in\mod 3), k \equiv 2 \,(\hskip -.1in\mod 6), 
\,\,8 \leq k \leq 20$, 
\vskip .1in
\hskip .7in $E_{3,\infty,k}^{-a}, \,\,\,\,1 \leq a \leq 97,
a \equiv 1\,  (\hskip -.1in\mod 3), L=2,\, k \equiv 2 \,(\hskip -.1in\mod 6),$
\newline \vskip 0in \hskip 1.3in $8 \leq k \leq 20$.
\vskip .1in
\hskip .7in $E_{3,\infty,k}^{-a}, \,\,\,\, 
1 \leq a \leq 98, k \equiv 4 \,(\hskip -.1in\mod 6),\,
4 \leq k \leq 22$.
\vskip .1in
The following functions deviate from rule (2):
\vskip .1in
\hskip .7in $	E_{2,\infty,k}^{-a}, \,\,\,\,1\leq a \leq 51, a \,\roman{odd}, 8 \leq k \leq 24,\, 
k \equiv 0\, (\hskip-.1in\mod 4)$. \vskip .1in
\flushpar Rule (2) predicts that  $\roman{ord}_2(c_{_0}[E_{2,\infty,k}^{-a}\,])
=3d_2(a)$ in this 
situation. Instead the constant terms 
obey the following rule :
$$\roman{ord}_2(c)=3d_2(a)+\roman{ord}_2(a+1)+k-5.\tag{3-1}$$

The following functions deviate from rule (3):
\vskip .1in
\hskip .7in $E_{3,\infty,k}^{-a}, \,\,\,\,1 \leq a \leq 49,
a \equiv 1 \,(\hskip -.1in\mod 3), k \equiv 0 \,(\hskip -.1in\mod 6), \,12 \leq k \leq 24$.
\vskip .1in
\flushpar  The weights of these functions are divisible by $3$, so
\,rule (3) predicts that $c_{_0}[E_{3,\infty,k}^{-a}\,] 
\equiv (-1)^a3^{d_3(a)} \,(\hskip -.1in\mod 3^{d_3(a)+1})$. Instead,
$$c_{_0}[E_{3,\infty,k}^{-a}\,] 
\equiv (-1)^{a+1}3^{d_3(a)} \,(\hskip -.1in\mod 3^{d_3(a)+1}).\tag{3-2}$$
The functions
\vskip .1in
\hskip .7in $E_{3,\infty,k}^{-a}, \,\,\,\,2 \leq a \leq 47,
a \equiv 2 \,(\hskip -.1in\mod 3), k \equiv 0 \,(\hskip -.1in\mod 6), \,12 \leq k \leq 24$
\vskip .1in
\flushpar also depart from rule (3).
In this situation, it is not true, as predicted by rule (3),
 that $\roman{ord}_3(c_{_0}[E_{3,\infty,k}^{-a}\,])=d_3(a)$. Instead
$$
 \roman{ord}_3(c_{_0}[E_{3,\infty,k}^{-a}\,])=d_3(a)+\roman{ord}_3(a+1)=
\delta \roman{\,\,(say)},\tag{3-3}
$$ 
We have not yet understood how these functions choose between the congruences
$c_{_0}[E_{3,\infty,k}^{-a}\,] 
\equiv \pm 3^{\delta} \,(\hskip -.1in\mod 3^{\delta+1})$, 
except that our data indicate that it depends only on the value of $a$.

The last set of functions in this survey deviating from rule (3) is:
\vskip .1in
\hskip .7in $E_{3,\infty,k}^{-a}, \,\,\,\,1 \leq a \leq 94,
a \equiv 1\,  (\hskip -.1in\mod 3), 
L=1,\, k \equiv 2 \,(\hskip -.1in\mod 6),$ \,\,\newline \vskip 0in
\hskip 1.3in $8 \leq k \leq 20$.
\vskip .1in \flushpar Here $w \equiv 1 (\hskip -.1in\mod 3)$, so rule (3) predicts that
$c_{_0}[E_{3,\infty,k}^{-a}\,] 
\equiv  3^{d_3(a)} \,(\hskip -.1in\mod 3^{d_3(a)+1})$. Actually for this set 
$$c_{_0}[E_{3,\infty,k}^{-a}\,] 
\equiv  -3^{d_3(a)} \,(\hskip -.1in\mod 3^{d_3(a)+1}).\tag{3-4}$$
\subhead{\rm 3.3} Divisibility properties of the Fourier coefficients of $j, \Delta$
 and their \vskip 0in \hskip .1in
reciprocals
\endsubhead

We observed a pattern of connections 
between corresponding Fourier coefficients 
(not the constant terms) of $1/\Delta$ and $j$, and between corresponding Fourier coefficients 
(not the constant terms) of $\Delta$ and $1/j$. These experiments
were motivated by the following considerations. Membership of $f^s$ in $C_2$ or $C_3$
for integers $s$, $1 \leq s \leq  B$ for some bound $B$ imposes conditions modulo powers 
of $2$ or $3$ 
on the Fourier coefficients of $f$ with exponent  $ \leq B-1$. It is easy to check,
for example, that if $f$ has a simple pole at infinity and $f^s \in C_2$ for $1 \leq s \leq 4$,
then
$$c_0[f] \equiv 8 \,(\hskip -.1in \mod 16),$$
$$c_1[f] \equiv 4 \,(\hskip -.1in\mod 8),$$
$$c_2[f] \equiv 0 \,(\hskip -.1in\mod 128),$$
$$c_3[f] \equiv 2 \,(\hskip -.1in\mod 4).$$
These calculations can be extended indefinitely. They suggest that there
is a systematic relationship between the 2- and 3-orders of
corresponding coefficients of any two functions satisfying the above requirements on $f$.
This led us to compare these orders in the expansions of $j$ and $1/\Delta$.
Let us denote $\roman{ord}_p(c_n[j]) - \roman{ord}_p(c_n[1/\Delta])$
as $\delta_{p,n}$.
For $-1 \leq n \leq 2470 \,\,(n \neq 0)$, we found that
$$n \equiv 0\,(\hskip -.1in \mod 2) \Rightarrow \delta_{2,n} = 3 \roman{ord}_2(n) + 1,$$
$$n \equiv 0\,(\hskip -.1in \mod 3) \Rightarrow \delta_{3,n} = 2 \roman{ord}_3(n) ,$$
$$n \equiv 1\,(\hskip -.1in \mod 3) \Rightarrow \delta_{3,n} = -1 ,$$
It is interesting to compare these rules with the congruences of Lehner 
([Lehner 1949], or [Apostol 1989], p.91). Writing $c(n)$ where we write $c_n[j]$,
they are:

$$ c(2^{\alpha}n) \equiv 0 \,(\hskip -.1in \mod 2^{3 \alpha + 8}),$$
$$ c(3^{\alpha}n) \equiv 0 \,(\hskip -.1in \mod 3^{2 \alpha + 3}),$$
$$ c(5^{\alpha}n) \equiv 0 \,(\hskip -.1in \mod 5^{\alpha + 1}),$$
$$ c(7^{\alpha}n) \equiv 0 \,(\hskip -.1in \mod 7^{\alpha}).$$
There are tantalizing hints of similar relations. For example, if  
$5 \leq n \leq 2470,$ \newline $ n \equiv 0 \,(\hskip -.1in  \mod 5), 
n \neq 2245,$ then $\delta_{5,n}=\roman{ord}_5(n),$ but 
 $\delta_{5,2245}=2.$

On general principles, we also compared $\Delta$ with $1/j,$  and found that
$\roman{ord}_p(c_n[1/j]) $ \newline $=\roman{ord}_p(c_n[\Delta]) \,\,$ 
for $ p=2,3, 1 \leq n \leq 4096.$ 
This also holds for $p=5$, 
if $n \not \equiv 3 $ or $4$ \newline
$ \,(\hskip -.1in \mod 5)$ and $1 \leq n \leq 1225$. 
These observations have some independent interest, because
 $c_n[\Delta]=\tau(n)$ (Ramanujan tau function). For example, one may imagine a proof of
Lehmer's conjecture that the Ramanujan tau function is non-vanishing,
 consisting of two parts: (1)
a proof of the above relations  for all positive $n$, 
and (2) a proof that $c_n[1/j]$ 
is non-vanishing. 

\vskip.1in
\head 4. Congruences \endhead
 The following scenario plays out only when we are lucky.
Given the power series of a modular form 
$f(x) = 1 + \sum_{n=1}^{\infty}a(n)x^n$, one uses M\"obius inversion 
and Apostol's Theorem 14.8 
to find the first few factors in the product expansion. 
 One then guesses 
the whole product expansion.  The product expansion then is 
used to guess how to write the
form as a monomial in Dedekind's $\eta$ function, 
and this relation is proved with the analytic theory
of modular forms. Then one derives the product expansion from that of $\eta$, 
and the 
recursion
among the Fourier coefficients using Apostol's Theorem 14.8. 
Finally, the recursion is used to
prove a special case of rules (1) - (3).

To illustrate, we will prove the 
following theorem, which is an example of rule (2):
\proclaim{Theorem 4.1} If $s =2^{\,x}, x=0,1,2,...,$ then 
$\,\,\roman{ord}_2\left(c_{\,0}\left[E_{\infty,4}^{-s}\right]\right)=3.$\endproclaim

We have written a similar proof for \proclaim{Theorem 4.2} 
If $s =  2^{\,x}D, D=1,3, or\, 5, x=0,1,2,...,$ then \,
$\Delta^{-s}$ lies in $C_2$.\endproclaim
\flushpar
which a reader can reproduce by imitating part of the proof of Theorem 4.1. 
The proof of Theorem 4.2 is simpler, because there is no need to derive 
the product expansions, but, as $D$ increases, it becomes messy. It seems that this process 
can be continued, but we have no reason to believe that it will work for every odd $D$. 
We would be surprised if  similar verifications of rule (1)
 could not also be written 
for $\roman{ord}_3\left(c_{\,0}\left[\Delta^{-s}\right]\right)$.
\vskip .1in \flushpar
\it Proof of Theorem 4.1. \rm 
The Fourier series of $E_{\infty,4}$ is
 monic integral, and therefore so are those of its integral powers. Thus 
the terms $R(n-a)$ on the right side of (2-21) are 
integral. If $s=2^x$, then $0<n<s$ implies that $\roman{ord}_2(n)<x$ . So (2-21) implies
 that 
$R(n) \equiv 0 \,(\hskip -.1in \mod 16)$. Also, by (2-21):
$$
R(s)=8\sum_{a=1}^s \sigma_1^{\roman{alt}}(a)R(s-a). \tag{4-1}
$$  

\flushpar
All the terms in the sum on the right side of (4-1), 
except the one corresponding to $a=s$, are
 congruent to zero modulo 16.
Therefore, $R(s)\equiv 8\sigma_1^{\roman{alt}}(s)R(0) \equiv 8(2^x+2^{x-1}+...+2-1)\cdot 1$ 
$\equiv 8 (\hskip-.1in\mod 16).$ Thus, $\roman{ord}_2(R(s))=3.$ 
But $R(s)=c_{\,0}[E_{\infty,4}^{-s}].$  
$\hskip 5in \boxed{} $

For the project of improving the bound in Theorem 2.1 in
the case $h \equiv 2 \, (\hskip -.07in \mod 4)$, 
the non-vanishing of the constant terms of the 
Fourier expansions of the $T_{2,h}$
forms is the key to our approach. We state some partial results in this direction
for $T$-series of both levels.
The arguments follow the approach used above and appear in [Brent 1994], 
Chapter 5.
\proclaim{Theorem 4.3} If $h \equiv
8 \, (\hskip -.07in \mod 12)$ and $r(1,h) = 2^x, x \geq 1$, then
$c_0[T_h] \equiv 16$ \newline 
$(\hskip -.07in \mod 32)$. If $h \equiv
2 \, (\hskip -.07in \mod 12)$ and $r(1,h) = 2^x, x \geq 1$, then
$c_0[T_h] \equiv
8 \, (\hskip -.07in \mod 32)$.
If $h = 2^x - 6 > 0$, then $c_0[T_{2,h}] \equiv
8 \, (\hskip -.07in \mod 16)$.  If $h = 2^x - 4 > 0$, then $c_0[T_{2,h}] \equiv
16$ \newline $(\hskip -.07in \mod 32)$.

\endproclaim

\head 5. Applications to the theory of quadratic forms \endhead
\subhead{\rm 5.1} Quadratic forms and modular forms \endsubhead
We tell how certain quadratic forms give rise to level two modular forms. 
For even $v$, set $\bold{x} = \, ^t(x_1,\,...,\, x_v)$, so that
$\bold{x}$ is a column vector. Let $A$ be an $v$ by $v$ square symmetric matrix 
with integer entries, even entries on the diagonal, and positive eigenvalues.
 Then $Q_A(\bold{x}) = \, ^t\bold{x} A \bold{x}$ is a homogenous second degree
polynomial in the $x_i$. We refer to $Q_A$ as the 
even positive-definite quadratic form associated to $A$.
If $\bold{x} \in \bold{Z}^v$, then $Q_A(\bold{x})$ is a non-negative even number,
which is zero only if $\bold{x}$ is the zero vector. The level of $Q_A$
is the smallest positive integer $N$ such that $NA^{-1}$ also has integer entries
and even entries on the diagonal.
 Let $\#Q_A^{-1}(n)$ denote 
the cardinality of
the inverse image in $\bold{Z}^v$ of an integer $n$ under the quadratic form $Q_A$.
\proclaim{Proposition 5.1} Suppose that $Q_A$ is a level two quadratic form.
Then the function 
$\Theta_A: \goth{H} \rightarrow \bold{C}$ satisfying
$$\Theta_A(z) = \sum_{n=0}^{\infty} \#Q_A^{-1}(2n)q^n\tag{5-1}$$
lies in $M(2,\frac v2)$.\endproclaim

\flushpar \it Proof. \rm We use machinery from [Miyake 1989].
Let $\chi: \bold{Z} \rightarrow \bold{C}$ be a Dirichlet character mod $N$,
 and $\alpha \in \Gamma_0(N)$ be the matrix 
$$\alpha = \left(\matrix a & b\\c & d\endmatrix \right) .$$
By abuse of notation we also let $\chi$ denote the 
character $\chi:\Gamma_0(N) \rightarrow
\bold{C}$ which acts by the map 
$\alpha \mapsto \chi(d).$
We have the stroke operator
$f|_h$: $$(f|_h \,\alpha )(z) = (cz+d)^{-h}f(\alpha \, z) \hskip 1in
(z \in \goth{H}).$$
We write $M(h, \Gamma_0(N), \chi)$ for the 
vector space of functions $f$ holomorphic on $\goth{H}^*$ such that
$f|_h \alpha = \chi(\alpha)f$ for all $\alpha \in \Gamma_0(N)$.
Thus $M(h, \Gamma_0(2), \chi)$ and $M(2,h)$ coincide for trivial $\chi$. 
The space $M(h, \Gamma_0(N), \chi)$ is itself trivial if 
$\chi(-1) \neq (-1)^h$. ([Miyake 1989], Lemma 4.3.2, p.115).
Thus the only non-trivial space $M(h, \Gamma_0(2), \chi)$
is $M(2,h)$.

Let $(n | m)$ be the Kronecker symbol. 
Let $A^{-1} = (b_{ij})$. We put $$\psi_A(m) = ((-1)^{v/2} det A | m)$$
and 
$$\Delta_A = \sum_{1 \leq i,j \leq v} b_{ij}\partial^2 / \partial x_i \partial x_j.$$
A spherical function of degree $\nu$ with respect to $A$
is a complex homogenous polynomial
$P(x_1,...,x_v) = P(\bold{m})$ (say) of degree $\nu$ annihilated by $\Delta_A$.
For $z \in \goth{H}$, let
$$ \theta_{A,P}(z) = \sum_{\bold{m} \in \bold{Z}^v} P(\bold{m}) \exp \left (2 \pi i 
\frac {Q_A(\bold{m})}2 z \right ).
$$
Then $\theta_{A,P} \in M(\frac v2 + \nu,\Gamma_0(2), \psi_A)$. ([Miyake 1989],
(3), p.192). Evidently, 
$$\Theta_A = \theta_{A,1} \in M \left (\frac v2, \Gamma_0(2), \psi_A\right).$$
In particular, $ M(\frac v2, \Gamma_0(2), \psi_A)$ is non-trivial, so
it must be $M(2,\frac v2)$. $\hskip .1in \boxed{} $

Since $M(2,h)$ is non-trivial only for even $h$, it also follows that $4 | v$.
\subhead{\rm 5.2} Quadratic minima \endsubhead
In this section we apply Theorem 2.1 to the problem of quadratic minima.
It is possible to improve the result slightly by an application of Theorem 4.3
to the sparse family of weights $h \equiv 2 \, (\hskip -.07in \mod 4)$ mentioned there.
It would be substantially improved by a proof that the constant term
of $T_{2,h}$ is non-zero for all $h \equiv 2 \,(\hskip -.07in \mod 4)$,
since this would improve Theorem 2.1.

\proclaim{Theorem 5.1} If $Q$ is a level two  
even positive-definite quadratic form
in $v$ variables, $8 | v$, then $Q$ represents
a positive integer $2n \leq 2 + \frac v4$. If $v \equiv 4$
 $(\hskip -.07in \mod 8)$, then $Q$ represents a positive integer 
$2n \leq 2 + \frac v2$. 
\endproclaim

\flushpar \it Proof. \rm Let $A$ be the matrix associated to $Q$, so that $Q=Q_A$. 
Suppose $v = 8u$. Then $\Theta_A \in M(2,4u)$ by Proposition 5.1. By Theorem 2.1, 
$\#Q_A^{-1}(2n) \neq 0$ for some 
$n, \, 1 \leq n \leq r(2,4u) = 1 + u$. That is, $Q$ represents
an integer $2n \leq 2(1 + u) = 2 + \frac v4$. On the other hand, suppose
 $v=8u+4$. Then $\Theta_A \in M(2,4u+2)$,
and $\#Q_A^{-1}(2n) \neq 0$ for some $n, \, 1 \leq n \leq 2r(2,4u+2) = 2(1 + u)$. Thus
$Q$ represents
an integer $2n \leq 4 + 4u = 2 + \frac v2$.$\hskip .1in \boxed{} $

\head 6. Conclusion \endhead

We  
don't know how to frame
natural descriptions of the families obeying the rules (1)-(3) from section 3.
We will only remark that some of our experiments indicate that
the arithmetic of the constant terms comes from the modularity of the underlying 
functions, but not from the properties of formal power series as they  relate to
 Ramanujan's congruences for the Ramanujan $\tau$ function. 
At the suggestion of Glenn Stevens, we formed non-modular series
obeying the Ramanujan congruences and checked the constant terms of their negative powers
without turning up examples of rules (1)-(3). It seems to be
the modularity of $\Delta$, for example, but not in a direct way
its obedience to the Ramanujan congruences,
that causes it to obey rule (1). 

On the basis of the observations reported in section 3, we could make many narrow conjectures.
Several seem to be worth stating.

\proclaim{Conjecture 1}

(i) The constant terms of the $T_{2,h}$ follow rule (2).

(ii) The forms $\Delta^{-s}$ and $j^s$, $s$ a positive integer, follow rule (1).

(iii) The forms $\Delta^{-1}$ and $j$ satisfy the relations of section 3.3 for
all non-zero 

\hskip .3in integers $n \geq -1.$

\endproclaim  

\flushpar Part (i) of Conjecture 1 implies
\proclaim{Conjecture 2} 
Suppose $f \in M(2,h)$ with Fourier expansion at infinity
$$f(z) = \sum_{n=0}^{\infty} A_n q^n, \,\,A_0 \neq 0.$$
If $h \equiv 2 \, (\hskip -.07in 
\mod 4)$, then some $A_n \neq 0, 1 \leq n \leq 1 + r(2,h)$.

\endproclaim

\flushpar In turn, Conjecture 2 implies

\proclaim{Conjecture 3}
If $Q$ is a level two or weakly level one even positive-definite quadratic form
in $v$ variables, $v \equiv 4$
 $(\hskip -.07in \mod 8)$, then $Q$ represents a positive integer 
$2n \leq 3 + \frac v4$. 
\endproclaim

\refstyle{C}

\enddocument